\documentclass[10pt]{article}
\usepackage{amsfonts}
\usepackage{mathrsfs}
\usepackage{amsmath}
\usepackage{amssymb}
\usepackage{graphicx}
\graphicspath{{figure/}}
\usepackage{amsthm,latexsym}
\renewcommand{\paragraph}{\roman{paragraph}}
\def \n{{\mathbf n}}
\def \m{{\mathbf m}}

\newtheorem{theorem}{\scshape \mdseries  Theorem}[section]
\newtheorem{lemma}[theorem]{\scshape \mdseries  Lemma}
\newtheorem{coro}[theorem]{\scshape \mdseries  Corollary}

\begin{document}

\title{\sf The least eigenvalues of signless Laplacian of non-bipartite graphs with pendant vertices\thanks{Supported by National Natural Science Foundation of China (11071002),
Program for New Century Excellent Talents in University,
Key Project of Chinese Ministry of Education (210091),
Specialized Research Fund for the Doctoral Program of Higher Education (20103401110002),
Science and Technological Fund of Anhui Province for Outstanding Youth  (10040606Y33),
Scientific Research Fund for Fostering Distinguished Young Scholars of Anhui University(KJJQ1001),
 Academic Innovation Team of Anhui University Project (KJTD001B),
 Fund for Youth Scientific Research of Anhui University(KJQN1003). }}
\author{Yi-Zheng Fan\thanks{Corresponding
author.  E-mail address: fanyz@ahu.edu.cn (Y.-Z. Fan),
wangy@ahu.edu.cn(Y. Wang), gh6596706@sina.com (H. Guo).}, \ Yi Wang, Huan Guo\\
  {\small  \it School of Mathematical Sciences, Anhui University, Hefei 230601, P. R. China} \\
 }
\date{}
\maketitle

\noindent {\bf Abstract:} In this paper we determine the graph whose least eigenvalue of signless Laplacian 
attains the minimum or maximum among all connected non-bipartite graphs of fixed order  and given number of pendant
vertices.  
Thus we obtain a lower bound and an upper bound for the least eigenvalue of signless Laplacian  of a graph in terms of the number of pendent vertices.

\noindent {\bf MR Subject Classifications:} 05C50, 15A18


\section{Introduction}

Let $G=(V,E)$ be a simple undirected graph with vertex set
$V=V(G)=\{v_1,v_2,\dots,v_n\}$ and edge set
$E=E(G)$.
The {\it adjacency matrix} of $G$ is
defined as the matrix $A(G)=[a_{ij}]$ of order $n$, where
$a_{ij}=1$ if $v_{i}$ is adjacent to $v_{j}$, and $a_{ij}=0$
otherwise.
The {\it degree matrix} of $G$ is defined by $D(G)={\rm
diag}\{d(v_1),d(v_2),\dots,d(v_n)\}$, where $d(v_i)$ is the degree of the vertex $v_i$.
The matrix $Q(G)=D(G)+A(G)$ is called the {\it signless Laplacian matrix} (or
{\it $Q$-matrix}) of $G$. It is known that $Q$ is nonnegative,
symmetric and positive semidefinite.
So its eigenvalues are all
nonnegative real numbers and can be arranged as:
$q_1(G)\geq q_2(G)\geq \cdots \geq q_n(G)\geq 0.$
We simply call the eigenvalues of $Q(G)$ as the {\it $Q$-eigenvalues} of the graph $G$,
and refer the readers to \cite{cve1, cve2, cve3, cve4, cve5} for the survey on this topic.
The least $Q$-eigenvalue $q_n(G)$ is denoted by $q_{\min}(G)$, and
the eigenvectors corresponding to $q_{\min}(G)$ are called the {\it first $Q$-eigenvectors} of $G$.

If $G$ is connected, then $q_{\min}(G)=0$ if and only if $G$ is bipartite.
So, the connected non-bipartite graphs are considered here. The very
early work on the least $Q$-eigenvalue can be found in \cite{des},
where the author discuss the relationship between the least
$Q$-eigenvalue and the bipartiteness of graphs.
Cardoso et al. \cite{car} and Fan et at. \cite{fan} investigate the
least $Q$-eigenvalue of non-bipartite unicyclic graphs.
 Liu et al.
\cite{liu} give some bounds for the clique number and independence
number of graphs in terms of the least $Q$-eigenvalue.
Lima et al. \cite{lima} survey the known results and present some
new ones for the least $Q$-eigenvalue. Our research group
\cite{wang} investigate how the least $Q$-eigenvalue  of a graph changes
 by relocating a bipartite branch from one vertex to another vertex, and
 minimize the least $Q$-eigenvalue among the connected graphs of fixed order which contain a given non-bipartite graph as an induced subgraph.

A graph is called  {\it minimizing} (or {\it maximizing}) in a class of graphs if its least $Q$-eigenvalue attains the minimum (or maximum) among all graphs in the class.
 Denote by $\mathscr{G}_n^k$ the set of connected non-bipartite graphs of order $n$ with $k$ pendant vertices.
In this paper we determine the unique minimizing graph and the maximizing graph in $\mathscr{G}_n^k$,
and hence provide a lower bound and an upper bound for the least $Q$-eigenvalue
of a graph in terms of the number of pendent vertices.

\section{Preliminaries}

We first introduce some notations.
We use $C_n$, $P_n$, $K_n$ denote the cycle, the path, the complete graph all on $n$ vertices, respectively.
We also use $Pv_1v_2\cdots v_n$ to denote a path on vertices $v_1,v_2,\ldots, v_n$ with edges $v_iv_{i+1}$ for $i=1,2,\ldots,v_{n-1}$.
Let $G$ be a graph.
The graph $G$ is called {\it trivial} if it contains only one vertex; otherwise, it is called {\it nontrivial}.
The graph $G$ is called {\it unicyclic}, if it is connected and has the same number of vertices and edges (or $G$ contains exactly one cycle).
The {\it girth} of $G$ is the minimum of the lengths of all cycles in $G$.
A {\it pendant vertex} of $G$ is a vertex of degree $1$.
A path $Pv_0v_1\cdots v_{t-1}v_t$ in $G$ is called a {\it pendant path} if
$d(v_1)=d(v_2)=\cdots=d(v_{t-1})=2$ and $d(v_{t})=1$.
If $t=1$, then $v_0 v_1$ is a  pendant edge of $G$.
In particular, if $d(v_0)\geq 3$, we say $P$ is a {\it maximal pendant path}.

Let $x=(x_1, x_2, \dots, x_n)^T $ be a column vector in
$\mathbb{R}^n$, and let $G$ be a graph on vertices $V(G)=\{v_1,v_2, \dots, v_n\}$.
The vector $x$ can be viewed as a function defined on
 $V(G)$, that is,  any vertex $v_i$ is given by the value
$x_i=:x{(v_i)}$. Thus the quadratic form $x^TQx$ can be written as
$$x^TQx=\sum\limits_{uv\in E(G)}[x(u)+x(v)]^2. \eqno(2.1)$$
One can find that $q$ is a $Q$-eigenvalue of $G$
corresponding to an eigenvector $x$ if and only if $x \neq 0$ and
$$[q-d(v)]x(v) =\sum_{u \in N_G(v)}x(u), \hbox{~ for each~} v\in V(G),\eqno(2.2)$$
where $N_G(v)$ denotes the neighborhood of the vertex $v$.
In addition, for an arbitrary unit vector $x \in \mathbb{R}^n$,
$$q_{\min}(G)\leq x^TQ(G)x,\eqno(2.3)$$
with equality if and only if $x$ is a first $Q$-eigenvector of $G$.

Let $G_1$ and $G_2$ be two vertex-disjoint graphs, and let $v \in
G_1$, $u \in G_2$. The {\it coalescence} of $G_1$ and $G_2$ with respect to $v$ and $u$, denoted
by $G_1(v) \diamond G_2(u)$, is obtained from $G_1$, $G_2$ by
identifying $v$ with $u$ and forming a new vertex.
Let $G$ be a connected graph, and let $v$ be a cut vertex of $G$.
Then $G$ can be expressed in the form $G=H(v)\diamond F(v)$, where $H$ and $F$ are subgraphs of $G$ both containing $v$.
Here we call $H$ (or $F$) a {\it branch of $G$ with root $v$}.
With respect to a vector $x$ defined on $G$, the branch $H$ is called {\it zero} if $x(v)=0$ for all $v \in V(H)$; otherwise $H$ is called {\it nonzero}.

 Let $G=G_1(v_2)
\diamond G_2(u)$, $G^*=G_1(v_1) \diamond G_2(u)$, where $v_1$ and
$v_2$ are two distinct vertices of $G_1$ and $u$ is a vertex of
$G_2$. We say $G^*$ is obtained from $G$ by {\it relocating $G_2$ from
$v_2$ to $v_1$}. In \cite{wang} the authors give some properties of the
first $Q$-eigenvectors, and discuss how the least $Q$-eigenvalue of a graph
changes when relocating a bipartite branch from one vertex to
another vertex; see the following results.

\begin{lemma} {\em\cite{wang}} \label{branch}
Let $H$ be a bipartite branch of a connected graph $G$ with root
$u$. Let $x$ be a first $Q$-eigenvector of $G$.\\
 {\em (1)} If $x(u)=0$, then
$H$ is a zero branch of $G$ with respect to $x$. \\
{\em (2)} If
$x(u)\neq 0$, then $x(p) \neq 0$ for every vertex $p$ of $H$.
Furthermore, for every vertex $p$ of $H$, $x(p)x(u)$ is positive or
negative, depending on whether $p$ is or is not in the same part of
bipartite graph $H$ as $u$; consequently, $x(p)x(q)<0$ for each edge
$pq \in E(H)$.
\end{lemma}

\begin{lemma} {\em\cite{wang}} \label{tree}
Let $G$ be a connected non-bipartite graph, and let $x$ be a first
$Q$-eigenvector of $G$. Let $T$ be a tree with root $u$, which is a
nonzero branch with respect to $x$. Then $|x(q)|< |x(p)|$ whenever
$p,q$ are vertices of $T$ such that $q$ lies on the unique path from
$u$ to $p$.

\end{lemma}

\begin{lemma} {\em\cite{wang}} \label{perturb}
Let $G_1$ be a connected graph containing at least two vertices
$v_1, v_2$, let $G_2$ be a connected bipartite graph containing a
vertex $u$. Let $G=G_1(v_2)\cdot G_2(u)$ and $G^{\ast}=G_1(v_1)\cdot
G_2(u)$. If there exists a first  $Q$-eigenvector of $G$ such that
$|x{(v_1)}|\geq{|x{(v_2)}|}$, then,
$$ q_{\min}(G^{\ast})\leq{ q_{\min}(G)}$$
\noindent{with equality only if $|x{(v_1)}|=|x{(v_2)}|$ and
$d_{G_2(u)}x{(u)}=-{\sum\limits_{v\in N_{G_2}(u)}{x(v)}}$}.
\end{lemma}

\begin{lemma}  {\em\cite{wang}} \label{perpath}
Let $G_1$ be a connected non-bipartite graph containing two vertices
$v_1,v_2$, and let $P$ be a nontrivial path with $u$ as an end vertex. Let
$G=G_1(v_2) \diamond P(u)$ and let $G^*=G_1(v_1) \diamond P(u)$. If
there exists a first $Q$-eigenvector $x$ of $G$ such that
 $|x(v_1)|>|x(v_2)|$ or $|x(v_1)|=|x(v_2)|>0$, then $$q_{\min}(G^*) <q_{\min}(G).$$
\end{lemma}

\section{Minimizing the least $Q$-eigenvalue among all graphs in $\mathscr{G}_n^k$}
Let $\mathscr{U}_n^k(g)$ denote the set of unicyclic
graphs of order $n$ with odd girth $g$ and $k \ge 1$ pendant vertices.
Denote by $U_n^k(g;l; l_1, l_2, \ldots, l_k) \in \mathscr{U}_n^k(g)$ the graph of order $n$ obtained by
coalescing $P_l$ with a cycle $C_g$ by identifying one of its end
vertices with some vertex of $C_g$,  and also coalescing this $P_l$ with each
of paths $P_{l_i}$ ($i=1, 2, \ldots, k$) by identifying its other end vertex with one of the end vertices of $P_{l_i}$,
where $l \ge 1$, $l_i \ge 2$ for $i=1, 2, \ldots, k$, $g+l+\sum_{i=1}^k l_i=n+k+1.$
If $l_1=l_2=\cdots=l_k=2$, $U_n^k(g; l; l_1, l_2, \ldots, l_k)$ is simply denoted by $U_n^k(g)$; see Fig. 3.1.

In this section, we first show that $U_n^k(g)$ is the unique minimizing graph in $\mathscr{U}_n^k(g)$, and then
investigate some properties of the least $Q$-eigenvalue and the corresponding eigenvectors of $U_n^k(g)$.
By the eigenvalue interlacing property (see following Lemma \ref{interlace}), the problem of determining the minimizing graph in $\mathscr{G}_n^k$
can be transformed to that of determining the minimizing graph in $\mathscr{U}_n^k(g)$.

\begin{theorem} \label{umin}
Among all graphs in $\mathscr{U}_n^k(g)$, $U_n^k(g)$ is the unique minimizing graph.
\end{theorem}

{\it Proof:} Let $G$ be a minimizing graph in $\mathscr{U}_n^k(g)$, and let $C_g$ be the unique
cycle of $G$ on vertices $v_1, v_2, \ldots, v_g$. The graph $G$ can be considered as  one
obtained from $C_g$ by identifying each $v_i$ with one vertex of some tree $T_i$ of order $n_i$ for
each $i=1,2, \ldots, g,$ where $\sum_{i=1}^g n_i=n$. Note that some trees $T_i$ may be trivial, i.e. $n_i=1$.

Let $x$ be a unit first $Q$-eigenvector of $G$. First, there exist
at least one $i$, $1 \leq i \leq g$, such that $x(v_i)\neq 0$;
otherwise, by Lemma \ref{branch}(1), each $T_i$, $1 \leq i \leq g$, is a zero
branch of $G$ with respect to $x$, and it follows that $x$ is the
zero vector, which is a contradiction.

We also assert that each nontrivial tree $T_j$ is a nonzero
branch with respect to $x$. Otherwise, there exists a nontrivial
tree $T_j$ attached at $v_j$, $1 \leq j \leq g$, such that $x(v_j) =
0$. By Lemma \ref{perturb}, relocating the tree $T_j$ from $v_j$ to $v_i$ for
some $i$ for which $x(v_i)\neq 0$, we obtain a graph in
$\mathscr{U}_n^k(g)$ with smaller least $Q$-eigenvalue.

Next, we contend all maximal pendant paths locate at the same vertex.
Otherwise, there exist two maximal pendant path, say $P$ and $P'$,
attached at $p$ and $p'$, respectively. Without loss of generality,
assume $|x(p)| \geq |x(p')| >0$. Note that $d(p') \geq 3$ by the
definition of maximal pendant path. Then by Lemma \ref{perpath}, we will arrive at
a new graph still in $\mathscr{U}_n^k(g)$ but with smaller least $Q$-eigenvalue by
relocating $P'$ from $q$ to $p$.
So $G$ is obtained from $C_g$ by attaching one path at some vertex of $C_g$ if $k=1$ (i.e. $G=U_n^k(g; n-g; 2)$),
or $G=U_n^k(g; l; l_1, l_2,\ldots, l_k)$ if $k \ge 2$ for some positive integers $l \ge 1$ and $l_i \ge 2 \;(i=1,2,
\ldots,k)$ satisfying $g+l+\sum_{i=1}^k l_i=n+k+1$.

To complete the proof, we only need to consider the case of $k \ge 2$ and prove that $l_1=\ldots=l_k=2$. If
not, say $l_i \geq 3$. Denote by $P_{l_i}=Pu_1u_2 \cdots u_{l_i}$, where $u_1$ is the common end vertex of other $k-1$ maximal pendant paths,
$d(u_1) = k+1, d(u_2)=\cdots=d(u_{l_i-1})=2, d(u_{l_i})=1.$ By Lemma \ref{tree} and above discussion, $0<|x(u_1)|<
|x(u_{l_i-1})|$.  Relocating some $P_{l_j}$ other than $P_{l_i}$ from $u_1$
to $u_{l_i-1}$, by Lemma \ref{perpath} we would arrive at a new graph in $\mathscr{U}_n^k(g)$ with smaller least
$Q$-eigenvalue, a contradiction. \hfill $\blacksquare$

\begin{coro} \label{usimple}
The least $Q$-eigenvalue of $U_n^k(g)$ has multiplicity one.
\end{coro}

{\it Proof:}
Let  $C_g$ be the unique cycle of $U_n^k(g)$, and let $v$
be the (unique) vertex lying on  $C_g$ with degree greater than $2$.
From the proof of Theorem \ref{umin}, the value of $v$ given by any first $Q$-eigenvector of $U_n^k(g)$ is nonzero.
Assume to the contrary, $x$ and $y$ are two linear independent first $Q$-eigenvectors of $U_n^k(g)$.
There exists a nonzero linear combination of $x$ and $y$ such that its value at $v$ equals zero, which yields a contradiction. \hfill $\blacksquare$

\begin{center}
\vspace{2mm}
\includegraphics[scale=0.65]{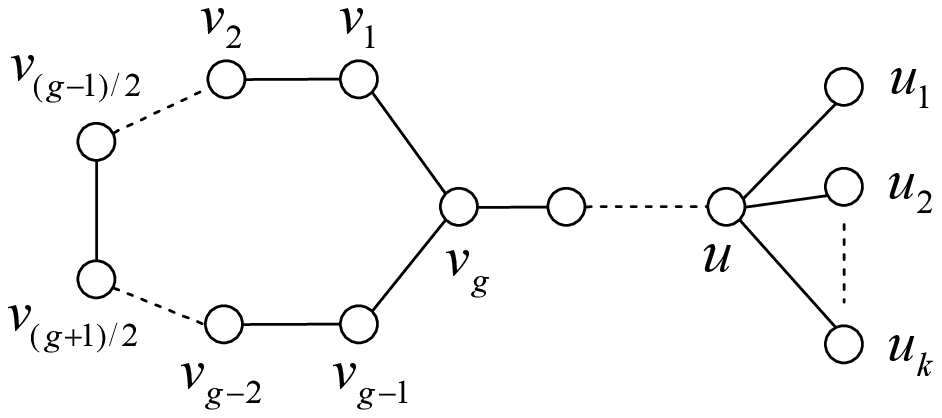}\\
\vspace{2mm}
 \small Fig. 3.1. The graph $U_n^k(g)$
\end{center}

\begin{lemma} \label{uvector}
Let $U_n^k(g)$ be the graph with some vertices labeled as in Fig. 3.1, where $v_1,v_2,\ldots,v_g$ are the vertices of
the unique cycle $C_g$ labeled in anticlockwise way. Let $x$ be a first $Q$-eigenvector of $U_n^k(g)$. Then\\
{\em (1)} $x(v_i)=x(v_{g-i})$ for $i=1,2, \ldots, \frac{g-1}{2}.$\\
{\em (2)} $x(v_{\frac{g-1}{2}})x(v_{\frac{g+1}{2}}) >0$, and $x(v_i)x(v_{i+1})<0$ for other edges $v_iv_{i+1}$ of $U_n^k(g)$ except $v_{\frac{g-1}{2}}v_{\frac{g+1}{2}}$.\\
{\em (3)} $|x(v_g)| > |x(v_1)|>|x(v_2)| >\ldots >
|x(v_{\frac{g-1}{2}})|>0.$
\end{lemma}

{\it Proof: } From the proof of Theorem \ref{umin}, the tree attached at $v_g$ is a nonzero branch with respect to $x$,
and by Lemma \ref{branch}(2) each edge $uv$ of the tree holds $x(u)x(v)<0$.
So it suffices to consider those edges on the cycle.

Observe that there exists an automorphism $\psi$ such that
$\psi(v_i)=\psi(v_{g-i})$ for $i=1, 2, \ldots, \frac{g-1}{2}$, and $\psi$
 preserves other vertices. Define a vector $x_{\psi}$ by $x_{\psi}(v)=x(\psi(v))$ for each
 vertex $v$ of $U_n^k(g)$. Then $x_{\psi}$ is also a unit first $Q$-eigenvector of $U_n^k(g)$.
Noting that $q_{\min}(U_n^k(g))$ is simple and $x_{\psi}(v_g)=x(v_g)\neq 0$, so
$x_{\psi}=x$, that is $x(v_i)=x(v_{g-i})$ for $i=1,2, \ldots,
\frac{g-1}{2}.$

Since $x(v_{\frac{g-1}{2}})=x(v_{\frac{g+1}{2}})$, we have
$x(v_{\frac{g-1}{2}})x(v_{\frac{g+1}{2}}) \geq 0$. If
$x(v_{\frac{g-1}{2}})=x(v_{\frac{g+1}{2}})=0$, by considering the
eigenvector equation (2.2) of $x$ at $v_{\frac{g-1}{2}}$, we have
$x(v_{\frac{g-3}{2}})=0$. Repeating the above discussion, we finally
obtain $x(v_g)=0$, a contradiction.

Next, we claim that $x(u)x(v) \leq 0$ for any edge $uv$ on the cycle $C_g$ other than $v_{\frac{g-1}{2}}v_{\frac{g+1}{2}}$.
Assume $pq$ is an edge on $C_g$ such that $x(p)x(q)>0$.
Partitioned the vertices of (the tree) $U_n^k(g)-v_{\frac{g-1}{2}}v_{\frac{g+1}{2}}$ into two parts 
 such that its edges join vertices from one part to vertices of the other part.
Note that $v_{\frac{g-1}{2}},v_{\frac{g+1}{2}}$ lie in the same part, and $p,q$ lie in different parts.
Define $\widetilde{x}$ on $U_n^k(g)-v_{\frac{g-1}{2}}v_{\frac{g+1}{2}}$ such that
$\widetilde{x}(v)=|x(v)|$ if $v \in S_1$, and $\widetilde{x}(v)=-|x(v)|$ if $v \in S_2$.
Then $\widetilde{x}^TQ(U_n^k(g))\widetilde{x}<x^TQ(U_n^k(g))x$, which yields a contradiction.
The remaining part of assertion (2) will be proved after showing the last assertion.

To prove the last assertion, we start with $|x(v_g)| > |x(v_1)|$.  If
not, relocating the pendant tree from $v_g$ to $v_1$, we can obtain a
graph $G'$ which holds $q_{\min}(G') \leq  q_{\min}(U_n^k(g))$ by Lemma \ref{perturb}.
Noting that $G'$ is isomorphic to $U_n^k(g)$, $q_{\min}(G') = q_{\min}(U_n^k(g))$.
 Also by lemma 2.3, the equality occurs only if $x(v_g)=-x(w)$, where $w$ is the neighbor of $v_g$ in the pendant tree.
 This contradicts the Lemma \ref{tree}.
 By induction,
assume that $|x(v_{i-1})|>|x(v_i)|$ for $1 \leq i \leq
\frac{g-1}{2}-1$, where $v_0:=v_g$.
By the eigenvector equation (2.2) of $x$ at $v_i$,
$$[q_{\min}(U_n^k(g))-2]x(v_i) =x(v_{i+1})+x(v_{i-1}).$$
Note that  $x(v_{i-1})x(v_i) \le 0, x(v_i)x(v_{i+1}) \leq 0$ by what we have proved, and
$0<q_{\min}(U_n^k(g))<1$ (see \cite{das}).
By the induction hypothesis, $|x(v_{i})|>|x(v_{i+1})|$, and the assertion (3) follows.
By the assertion (3) we now can deduce the assertion (2). \hfill $\blacksquare$

\begin{coro} \label{unonzero}
Let $x$ be a first $Q$-eigenvector of $U_n^k(g)$. Then $x$ contains
no zero entries.
\end{coro}

Denote by $\alpha_n^k(g)$ the minimum of the least $Q$-eigenvalues of graphs in $\mathscr{U}_n^k(g)$, that is, the least $Q$-eigenvalue of
$U_n^k(g)$.

\begin{lemma} \label{umonotone}
$\alpha_n^k(g)$ is strictly increasing with respect to $k \geq 1$ and  odd $ g \geq
3 $, respectively.
\end{lemma}

{\it Proof:}  Let $U_n^k(g)$ have some vertices labeled as in Fig. 3.1.
Let $x$ be a first $Q$-eigenvector of $U_n^k(g)$.
Suppose $k \geq 2$. Replacing the edge $uu_2$ by $u_1u_2$, we arrive at a new graph
$\overline{G}  \in \mathscr{U}_n^{k-1}(g)$, which holds that $q_{\min}(\overline{G}) < q_{\min}(U_n^k(g))$ by Lemma \ref{perturb} as $|x(u)| < |x(u_1)|$.
So, by Theorem \ref{umin} we have
 $$\alpha_n^{k-1}(g) \le q_{\min}(\overline{G}) < q_{\min}(U_n^k(g))=\alpha_n^k(g).$$

Next we prove the second result. Suppose $g \geq 5$.  Replacing the edge $v_{g-2}v_{g-1}$ by edge $v_{g-2}v_1$,
we obtain a new graph  $\widetilde{G} \in \mathscr{U}_n^{k+1}(g-2)$ whose least $Q$-eigenvalue is not greater than $\alpha_n^k(g)$ as $x(v_1)=x(v_{g-1})$ by Lemma \ref{uvector}.
 Then
 $$\alpha_n^k(g-2) < \alpha_n^{k+1}(g-2) \leq  q_{\min}(\widetilde{G}) \leq \alpha_n^k(g).$$ The result
 follows. \hfill $\blacksquare$

\begin{lemma} {\em\cite{car}} \label{interlace}
Let $G$ be a graph of order $n$ containing an edge $e$. Then
$$q_1(G-e) \leq q_1(G) \leq q_2(G-e) \leq q_2(G) \leq q_3(G-e)
\leq \ldots \leq q_n(G-e) \leq q_n(G).$$
\end{lemma}

Now  we arrive at the main result of this section.
\begin{theorem} \label{mineig}
Among all graphs in $\mathscr{G}_n^k$, $U_n^k(3)$ is the unique minimizing graph.
\end{theorem}

{\it Proof:} Let $G$ be a minimizing graph in $\mathscr{G}_n^k$. Then $G$
contains at least an induced odd cycle, say $C_g$.
Let $G'$ be a connected unicyclic spanning subgraph of $G$, which contains $C_g$ as the unique cycle and contains all pendant edges of $G$.
Thus $G' \in \mathscr{U}_n^{k'}(g)$, where $k' \ge k$.
By Lemma 3.6 and Lemma 3.5,
$$q_{\min}(U_n^k(3))=\alpha_n^k(3) \le \alpha_n^{k}(g) \leq \alpha_n^{k'}(g) \leq
q_{\min}(G') \leq q_{\min}(G).\eqno(3.1)$$
As $G$ is a minimizing graph in $\mathscr{G}_n^k$, all inequalities in (3.1) hold as equalities,
which implies $k'=k, g=3$ and $G'=U_n^k(3)$ by Lemma 3.5 and Theorem \ref{umin},  and also $q_{\min}(G)=q_{\min}(U_n^k(3))$.

Now we return to the origin graph $G$, which is obtained from $G'=U_n^k(3)$ possibly by adding some edges.
Suppose $E(G)\setminus E(U_n^k(3)) \neq \emptyset$.
Recalling the definition of $G'$ and $U_n^k(3)$, the set $E(G)\setminus E(U_n^k(3))$ consists of some edges
 joining the vertices of $C_3$ and the vertices of $P_l$ or some edges within the vertices of $P_l$.
So, for each edge $uv \in E(G)\setminus E(U_n^k(3))$, if $x$ is a first $Q$-eigenvector of $U_n^k(3)$, then $x(u)+x(v) \ne 0$ by Lemma \ref{uvector}(3) and Lemma \ref{tree}.

Let $x$ be a unit first $Q$-eigenvector of $G$. Then
\begin{eqnarray*}
q_{\min}(G)&=&\sum_{uv \in E(G)}[x(u)+x(v)]^2\\
&=&\sum_{uv \in E(U_n^k(3))}[x(u)+x(v)]^2 + \sum_{uv \in E(G)\setminus
E(U_n^k(3))}[x(u)+x(v)]^2\\
&\geq&\sum_{uv \in E(U_n^k(3))}[x(u)+x(v)]^2 \geq q_{\min}(U_n^k(3))
\end{eqnarray*}
Since $q_{\min}(G)=q_{\min}(G')$, $x$ is also an first
$Q$-eigenvector of $U_n^k(3)$, and for each edge $uv \in E(G)\setminus E(U_n^k(3))$, $x(u)+x(v) = 0$, which yields a contradiction. The result
follows.  \hfill $\blacksquare$

By Theorem \ref{mineig} and Lemma \ref{umonotone}, we have the following result.

\begin{coro} \label{eigpen}
Let $G$ be a connected graph of order $n$ which contains pendant vertices.
Then $q_{\min}(G) \ge q_{\min}(U_n^1(3))$
with equality if and only if $G=U_n^1(3)$.
If, in addition, $G$ contains $k$ pendant vertices,
then
$q_{\min}(G) \ge q_{\min}(U_n^k(3))$
with equality if and only if $G=U_n^k(3)$.
\end{coro}

Cardoso et al. \cite{car} determine the unique minimizing graph of non-bipartite connected graph, i.e. the graph $U_n^1(3)$.
By Lemma \ref{interlace}, the minimizing graph is unicyclic.
If we know that the (unicyclic)  minimizing graph contains pendant vertices, then
we also determine this minimizing graph by Corollary \ref{eigpen}.

\section{Maximizing the least $Q$-eigenvalue among all graphs in $\mathscr{G}_n^k$}
Let $\n=(\n_1,\n_2, \ldots, \n_{n-k}) \in \mathbb{N}^{n-k}$ be a nonnegative integer sequence arranged in
non-increasing order, where $\n_1+\n_2+ \cdots+ \n_{n-k}=k$.
In this section, all nonnegative integer sequence has the same form as $\n$.
Denote by $K(\n)$ the graph obtained from $K_{n-k}$ on vertices $v_1,v_2,\ldots,v_{n-k}$ by
attaching $\n_i$ pendant edges to $v_i$ for $i=1,2, \ldots n-k$, respectively.
By Lemma 3.6, the maximizing graph in $\mathscr{G}_n^k$ is achieved
by $K(\n)$ for some $\n \in \mathbb{N}^{n-k}$.

\begin{lemma} \label{maxvec}
Let $x$ be a first $Q$-eigenvector of $K(\n)$. If $\n_i > \n_j$, then
$|x(v_i)| \ge |x(v_j)|$.
\end{lemma}

{\it Proof:} Assume to the contrary, $|x(v_i)| < |x(v_j)|$.
Relocating $\n_i-\n_j$ pendent edges from $v_i$ to $v_j$,
by Lemma \ref{perturb}, $q_{\min}(\widetilde{G}) < q_{\min}(K(\n))$.
But $ \widetilde{G}$ is isomorphic to $K(\n)$ so that $q_{\min}(\widetilde{G}) = q_{\min}(K(\n))$, a contradiction.\hfill $\blacksquare$

 Recalling the notation of majorization, if
$\n=(\n_1, \n_2, \ldots, \n_r)$ and $\m=(\m_1,\m_2,
\ldots,\m_r)$ are two nonnegative integer sequences arranged in
non-increasing order, then $\n$ {\it majorizes} $\n$, denote by
$\n \succeq \m$, if, for $1 \le k \le r-1$, $$\sum_{i=1}^k\n_i \geq \sum_{i=1}^k \m_i {\rm~~and~~}  \sum_{i=1}^r\n_i=\sum_{i=1}^r\m_i.$$
If $\n \succeq \m$ and $\n \ne \m$, we will denote $\n \succneqq \m$.

\begin{lemma} \label{maj}
Let $\n=(\n_1, \n_2, \ldots, \n_{n-k})$ be a nonnegative integer sequences arranged in
non-increasing order, where $\n_1+\n_2+ \cdots+ \n_{n-k}=k$.
If $\n_1-\n_{n-k} \ge 2$, there exists a nonnegative integer sequences $\m \in \mathbb{N}^{n-k}$
such that
$ \n \succneqq \m$ and $q_{\min}(K(\n)) \le q_{\min}(K(\m)).$
\end{lemma}

{\it Proof:}
Suppose $K(\n)$ has the vertices labeled at the beginning of this section.
Relocating a pendant edge from $v_1$ to $v_{n-k}$, we will arrive at a new graph $G$ isomorphic to $K(\m)$ for some $\m\in \mathbb{N}^{n-k}$.
Surely $\n \succneqq \m$.
Let $x$ be a first $Q$-eigenvector of $G$.
By Lemma \ref{maxvec}, $|x_{v_1}| \ge |x_{v_{n-k}}|$ as $\n_1-1 \ge \n_{n-k}+1$.
(If $\n_1-1 = \n_{n-k}+1$ and $|x_{v_1}| < |x_{v_{n-k}}|$, we may interchange the labeling of $v_1,v_{n-k}$ to make the above inequality hold.)
Now relocating a pendant edge from $v_{n-k}$ to $v_1$, we go back to the original graph $K(\n)$.
By Lemma \ref{perturb}, $q_{\min}(K(\n)) \le q_{\min}(G)=q_{\min}(K(\m))$. The result holds.\hfill $\blacksquare$

By repeatedly using Lemma \ref{maj}, we get the following result.

\begin{theorem}\label{max}
The maximizing graph in $\mathscr{G}_n^k$ can be achieved by $K(\n)$, where
$\n=(\n_1,\n_2, \ldots, \n_{n-k})$, $\sum_{i=1}^{n-k}\n_i=k$, and $|\n_i-\n_j|\leq 1$ for all $i,j=1,2,\ldots,n-k$.
\end{theorem}

\begin{coro} \label{upb}
Let $G$ be a connected graph containing pendant vertices.
Then
$$q_{\min}(G) \le \frac{n-1+\frac{1}{n-1}-\sqrt{n^2-6n+11+\frac{1}{(n-1)^2}}}{2}.$$
If, in addition, $G$ contains $k \ge 1$ pendant vertices,
then
$$q_{\min}(G) \le \frac{n-k+\frac{k}{n-k}-\sqrt{(n-k-2)^2+2k+\frac{k^2}{(n-k)^2}}}{2}.$$
\end{coro}

{\it Proof:}
Assume $G \in \mathscr{G}_n^k$ for some $k \ge 1$.
By Theorem \ref{max}, $q_{\min}(G) \le q_{\min}(K(\n))$, where $\n$ has the prescribed property in Theorem \ref{max}.
Let $t:=\lceil k/(n-k) \rceil$, and let $B$ be the principal submatrix of $Q(K(\n))$ indexed by the vertex with degree $t+n-k-1$ and $t$ pendant vertices adjacent to it.
By the eigenvalue interlacing property of symmetry matrices,
\begin{align*}
q_{\min}(K(\n)) \le q_{\min}(B)&=\frac{n-k+t-\sqrt{(n-k+t)^2-4(n-k-1)}}{2}\\
&\le \frac{n-k+\frac{k}{n-k}-\sqrt{(n-k+\frac{k}{n-k})^2-4(n-k-1)}}{2}\\
&=\frac{n-k+\frac{k}{n-k}-\sqrt{(n-k-2)^2+2k+\frac{k^2}{(n-k)^2}}}{2},
\end{align*}
where  $q_{\min}(B)$ is the least eigenvalue of $B$.
Noting the function $f(k):=\frac{n-k+\frac{k}{n-k}-\sqrt{(n-k-2)^2+2k+\frac{k^2}{(n-k)^2}}}{2}$ is strictly decreasing with respect to $k$,
so we get the first result of this theorem. \hfill $\blacksquare$

We now give a remark on Lemma \ref{maxvec}, Lemma \ref{maj} and Theorem \ref{max}.
Consider the graphs $K(2,2,2,0)$ and $K(2,2,2,1)$ in Fig. 4.1.
Using the software {\sc Mathematica}, we find they have the same least $Q$-eigenvalues, both being $(5-\sqrt{17})/2$ with multiplicity $2$.
So, the inequality in Lemma \ref{maj} may hold as an equality; and
the maximizing graph in $\mathscr{G}_n^k$ may not be unique.
The two independent first $Q$-eigenvectors of $K(2,2,2,0)$ are listed below:
{\small $$x=\left(\frac{\sqrt{17}-3}{2}, 0, \frac{3-\sqrt{17}}{2}, 0, -1, -1, 0, 0, 1,1\right), \ \
y=\left(\frac{\sqrt{17}-3}{2}, \frac{3-\sqrt{17}}{2}, 0,0, -1, -1, 1, 1, 0,0\right).$$}
We find that $|x(v_1)|>|x(v_2)|=|x(v_4)|$ even though $\n_1=\n_2>\n_4$.
So, in Lemma \ref{maxvec}, if $\n_i > \n_j$ we cannot say $|x(v_i)| > |x(v_j)|$,
and if $\n_i=\n_j$ we also cannot say $|x(v_i)| = |x(v_j)|$.

\begin{center}
\vspace{2mm}
\includegraphics[scale=0.65]{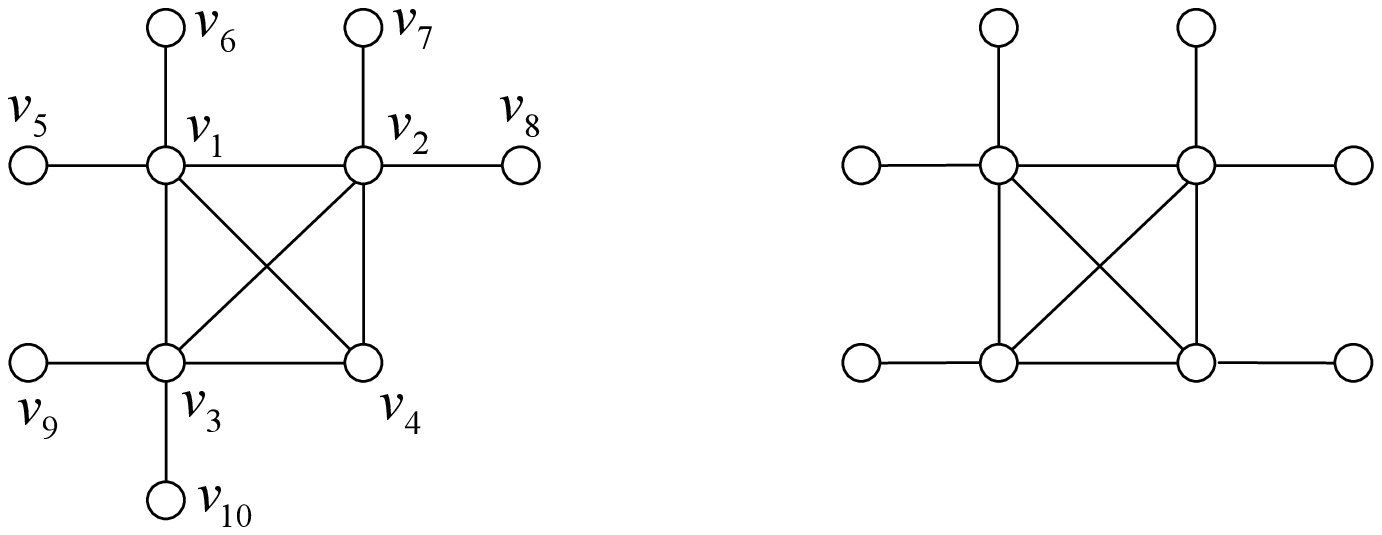}\\
\vspace{2mm}
 \small Fig. 4.1 The graph $K(2,2,2,0)$ (left) and $K(2,2,1,1)$ (right)
\end{center}

Finally we give a remark on some upper bounds of the least $Q$-eigenvalue of a graph $G$ in terms of minimum degree $\delta(G)$.
Liu and Liu \cite{liul} observe that $q_{\min}(G) \le \delta(G)$.
Das \cite{das} show that $q_{\min}(G) < \delta(G)$.
Lima et al. \cite{lima} improve the bound as
$$q_{\min}(G) \le \frac{n-1+\delta(G)-\sqrt{(n-1-\delta(G))^2+4}}{2}<\delta(G).$$
If the graph $G$ contains pendant vertices, i.e. $\delta(G)=1$, then
the above bound is
$$\frac{n-\sqrt{n^2-4n+8}}{2} > \frac{n-1+\frac{1}{n-1}-\sqrt{n^2-6n+11+\frac{1}{(n-1)^2}}}{2}.$$
So we give a subtle upper bound for the least  $Q$-eigenvalue of a graph if the graph contains pendant vertices.

\end{document}